\documentclass[12pt]{amsart}
\usepackage[arrow,dvips,matrix,arrow,ps,color,line,curve,frame]{xy}

\usepackage{amssymb}

\advance\textheight1cm

\let\oldlabel=\label
\def\prellabel{\marginparsep=1em
    \def\label##1{\oldlabel{##1}\ifmmode\else\ifinner\else
         \marginpar{{\footnotesize\ \\ \tt
                    ##1}}\fi\fi}}
\advance\headheight1.15pt
\let\:\colon
\let\epsilon\varepsilon

\let\phi=\varphi
\let\theta=\vartheta

\let\Bbb=\mathbb

\def\opn#1#2{\def#1{\operatorname{#2}}}
\opn\gp{gp} \opn\Max{Max} \opn\Ker{Ker} \opn\Coker{Coker}
\opn\Ext{Ext} \opn\conv{conv} \opn\chara{char} \opn\n{n} \opn\h{h}
\opn\GL{GuL} \opn\SL{SL} \opn\sn{sn} \opn\inte{int} \opn\End{End}
\opn\rank{rank} \opn\Aff{Aff} \opn\Spec{Spec} \opn\Proj{Proj}
\opn\gr{gr} \opn\QF{QF} \opn\I{Im} \opn\Hom{Hom} \opn\Aut{Aut}

\opn\calg{Alg^c} \opn\m{Mod}  \opn\pr{pr} \opn\E{E}
 \def\gac{GA^c} \opn\nlc{{\bf NL}^c}
\def\gec{GE^c} \def\gic{GI^c} \opn\U{U} \opn\X{\overline X}
\opn\Y{\overline Y} \opn\B{B} \opn\Mod{Mod} \opn\i{I/I^2}

\def\m{{\mathfrak M}}
\opn\w{w} \opn\inte{int} \opn\pyr{pyr} \opn\l{l} \opn\r{r}
\opn\const{const}

\def\ZZ{{\Bbb Z}}

\def\NN{{\Bbb N}}
\def\QQ{{\Bbb Q}}

\def\Q{{\Box\kern1pt}}%

\def\PP{{\Bbb P}}

\def\1{^{-1}}

\newtheorem{lemma}{Lemma}[section]
\newtheorem{corollary}[lemma]{Corollary}
\newtheorem{theorem}[lemma]{Theorem}

\theoremstyle{definition}

\newtheorem{remark}[lemma]{Remark}

\begin{document}

\title[commutators automorphisms]
{Commutator automorphisms of\\ formal power series rings}

\author[Joseph Gubeladze and Zaza Mushkudiani]
{Joseph Gubeladze$^{\diamond}$ and Zaza Mushkudiani}

\thanks{$^\diamond$ Corresponding author}
%\thanks{$^\bullet$ Supported by ***}


\address{Department of Mathematics, San Francisco
State University, San Francisco, CA 94132, USA}

\email{soso@math.sfsu.edu}

\address{TBC Bank, 11 Chavchavadze Ave., 0179 Tbilisi, Georgia}


\begin{abstract}
For a big class of commutative rings $R$ every continuous
$R$-auto\-mor\-phism of $R[[X_1,\ldots,X_n]]$ with the identity
linear part is in the commutator subgroup of
$\Aut(R[[X_1,\ldots,X_n]])$. An explicit bound for the number of
the involved commutators and a $K$-theoretic interpretation of
this result are provided.
\end{abstract}

\maketitle

\section{Introduction}\label{statement}

\subsection{Motivation}
A $K$-theoretical framework for some well known conjectures on
automorphisms and idempotent endomorphisms of polynomial rings was
initiated at the end of the 1970s in the works \cite{BW,C1,C2,CW}.
However, the development since then has been blocked due to
intractability of these conjectures. A recent evidence of this
intractability is provided by \cite{SU}.

There are two natural variations of the mentioned nonlinear
$K$-theory that makes things easier and tractable. One variant
corresponds to extending polynomial algebras to \emph{polytopal
algebras} -- essentially affine cones over projective toric
varieties -- and restricting homomorphisms to graded
homomorphisms. A systematic theory in this direction, including
higher groups, is developed in \cite{BrG1,BrG2,BrG3,BrG4}.

Here we consider the other variation that corresponds to the
completion process and leads to the category of formal power
series and their continuous homomorphisms. For higher groups this
approach leads to a challenging conjectural homological
computation -- the isomorphisms (\ref{conj}) below for which
supporting low dimensional results are obtained in this paper.

The freeness of projective objects in the context of complete
algebras was proved in \cite{T}, as an application of a technique
designed to study retracts of polynomial algebras, and also in
\cite{P}, in the general context of projective algebras over
nilpotent algebraic theories. One should also mention the work
\cite{K}, establishing an analogue of the $K$-theoretical
local-global principle for complete algebras.

Below we show (Theorem \ref{main}) how the standard method of
`indeterminate coefficients' implies the aforementioned result on
retracts of $R[[X_1,\ldots,X_n]]$. Furthermore, the same method
almost immediately yields the corresponding $K_1$-analogue for the
group $\Aut(R[[X_1,\ldots,X_n]])$ if the ground ring $R$ contains
$\QQ$ (Theorem \ref{Main}).

Main emphasize in this paper is put on derivation of a similar
result on the automorphism group when the coefficient ring is an
algebra over a finite field (Theorem \ref{MAIN}). Here the things
are a bit complicated. One needs to invoke a complex induction on
coefficient indices, depending on the residues modulo the
characteristic. The main difficulty to make the indeterminate
coefficients method work has been guessing many auxiliary
automorphisms yielding feasible systems of equation. Also, one
needs to treat separately the univariate case and the passage to
the multivariate case. For a variant of the result in the case
$n=1$ and $R=\ZZ_p$, where $p$ is a prime number $\not=2,3$, see
\cite{Ca1}. The corresponding groups of automorphisms are related
to \emph{Nottingham groups} -- objects of recent active
investigation \cite{Ca2}.

We have not been able to involve algebras over $\ZZ_2$ and $\ZZ_3$
in our consideration.

The higher $K$-theory machinery of symmetric monoidal categories
\cite{Gr} explains what the direct higher analogue of the main
result (after interpreting it as a claim on the fundamental group
of certain $K$-theoretic spaces) should be. We conjecture the
isomorphisms of integral homologies
\begin{equation}\label{conj}
H_i(E(R),\ZZ)=H_i(\left[\gac(R),\gac(R)\right],\ZZ),\quad i\ge2
\end{equation}
where $E(R)$ is the stable group of elementary matrices and
$\gac(R)$ stands for the formal non-linear analogue of the general
linear group $GL(R)$. These isomorphisms seem to be highly
nontrivial already for $i=2$. In this particular case a theory of
universal Steinberg relations should be possible, generalizing
Milnor's $K_2$. In the polytopal $K$-theory this is done in
\cite{BrG3}.

\subsection{Main results} Throughout the paper $R$ is a commutative ring.

$\Mod(R)$ denotes the category of $R$-modules, and $\PP(R)$
denotes the full subcategory of finitely generated projective
modules. Let $\calg(R)$ be the category of augmented $R$-algebras,
complete in the additive topology of the augmentation ideal.   For
$A\in\calg(R)$ let $I(A)=\Ker(A\to R)$. The coproduct of two
algebras in $\calg(R)$ is given by the additive completion of the
tensor product: $A\sqcup B=(A\otimes_RB)^\wedge_{I(A\otimes B)}$,
$I(A\otimes B)=\I\left((I(A)\otimes B)\ \oplus\ (A\otimes I(B))\to
A\otimes B\right)$. This makes $\calg(A)$ a symmetric monoidal
category. We have the functors:
\begin{align*}
\hat S:\Mod(R)\to\calg(R),\quad \hat S(M)&=(R\oplus M\oplus S^2(M)
\oplus\cdots)^\wedge_{0\oplus M\oplus S^2(M)
\oplus\cdots}\\
\i:\calg(R)\to\Mod(R),\quad \i(A)&=I(A)/I^2(A).
\end{align*}

 Let $\X=(X_1,\ldots,X_n)$. A \emph{retract} of $R[[\X]]$ is an algebra $A\in\calg(R)$ for
which there are morphisms $A\to R[[\X]]\to A$ whose composite is
${\bf 1}_A$.

\begin{theorem}\label{main}
Any retract of $R[[\X]]$ is of type $\hat S(P)$ for some
$P\in\PP(R)$.
\end{theorem}

Denote by $\gac_n(R)$ the group of augmented continuous
$R$-automorphisms of the formal power series ring $R[[\X]]$, where
`augmented' and `continuous' are understood with respect to the
ideal $(\X)R[[\X]]\subset R[[\X]]$ and the corresponding additive
topology on $R[[\X]]$. Such an automorphism is determined by its
values on the variables
$$
X_i\mapsto h_i^{(1)}+h_i^{(2)}+\cdots,\quad i=1,\ldots,n
$$
where $h_i^{(j)}\in R[\X]$ is a homogeneous polynomial of degree
$j$ for $i=1,\ldots,n$.

The surjective group homomorphism
$\gac_n(R)\to\Aut(R[[\X]]/\m^2)=GL_n(R)$ is split by the natural
embedding $GL_n(R)\to\gac_n(R)$ whose image consists of the linear
transformations of the $X_i$. Let $\gic_n(R)$ and $\gec_n(R)$
denote respectively the preimages of the trivial and elementary
subgroups: $\{{\bf 1}\}\subset E_n(R)\subset GL_n(R)$.

For a natural number $k$ and a group $\Gamma$ we let
$[\Gamma,\Gamma]^{(k)}$ denote the subset
$$
\{[\beta_1,\gamma_1]\cdots[\beta_k,\gamma_k]\ :\
\beta_1,\gamma_1,\ldots,\beta_k,\gamma_k\in\Gamma\}\subset
[\Gamma,\Gamma]
$$
where $[-,-]$ refers to the commutator.

\begin{theorem}\label{Main}
If there exist $b,c\in\U(R)$ such that $(b^m-1)R+(c^m-1)R=R$ for
all natural numbers $m$ then
$\gic_n(R)\subset[\gac_n(R),\gac_n(R)]^{(1)}$.
\end{theorem}

The class of rings in Theorem \ref{Main} contains all
$\QQ$-algebras.

Let $p$ be a prime number $\ge5$.

\begin{theorem}\label{MAIN}
(a) If $\ZZ_p\subset R$ then
$\gic_1(R)=[\gac_1(R),\gac_1(R)]^{(2)}$.

(b) If $\ZZ_p\subset R$ and $n\ge2$ then
$\gic_n(R)\subset[\gac_n(R),\gac_n(R)]^{(n(n+3))}$.

(c) If $R$ is a field of characteristic $p$ then
$\gic_n(R)\subset[\gac_n(R),\gac_n(R)]^{(4n)}$.
\end{theorem}

\

\noindent\emph{Convention}. Below $R^n$ is thought of as the
module of $n$-columns, $\U(R)$ refers to the group of units of
$R$. For a univariate formal power series $f(X)\in R[[X]]$ its
$m$-th coefficient is denoted by $f(X)_m$.

\section{Proof of Theorem \ref{main}}\label{Retracts}
Let $A\subset R[[\X]]$ be an $R$-subalgebra for which there exists
an $R$-algebra homomorphism $\pi:R[[\X]]\to A$ such that
$\pi|_A={\bf 1}_A$. We put
$$
\pi(X_i)=f_i=f_i^{(1)}+f_i^{(2)}+\cdots,\quad i=1,\ldots,n
$$
where $f_i^{(k)}\in R[\X]$ is a homogeneous polynomial of degree
$k$. Let $M(\pi)$ denote the $n\times n$-matrix whose $i$th column
is $(a_{i1},\ldots,a_{in})^{\intercal}$ where
$f^{(1)}_i=a_{i1}X_1+\cdots+a_{in}X_n$. The submodule of $R^n$
generated by the columns of $M(\pi)$ is naturally identified with
the module $P=(I/I^2)(A)$. The latter is a projective $R$-module:
we have the split epimorphism $(I/I^2)(\pi):R^n\to P$. Put
$Q=\Ker\left((I/I^2)(\pi)\right)$ and fix an epimorphism
$\rho:R^n\to Q$ split by the embedding $Q\subset R^n$. Then we
have the following split $R$-algebra epimorphism
$$
\pi\sqcup\hat S(\rho):R[[\X,\Y]]\to A\sqcup\hat S(Q).
$$
(Here $\Y=(Y_1,\ldots,Y_n)$.) Because $P\oplus Q\simeq R^n$ there
exists $\alpha\in\Aut_R(R^{2n})=GL_{2n}(R)$ such that
$\alpha(P\oplus Q)=R^n\oplus0$. Then the idempotent endomorphism
$\kappa=\hat S(\alpha)\left(\pi\sqcup\hat S(\rho)\right)\hat
S(\alpha)^{-1}:R[[\X,\Y]]\to R[[\X,\Y]]$ is of the form
\begin{align*}
\kappa(X_i)=
X_i+g^{(2)}_i+g^{(3)}_i+\cdots,\quad i=1,\ldots,n,\\
\kappa(Y_j)=h^{(2)}_j+h^{(3)}_j+\cdots,\quad j=1,\ldots,n\
\end{align*}
where $g_i^{(k)},h_j^{(k)}\in R[\X,\Y]$ are homogeneous
polynomials of degree $k$. We want to show the equality
\begin{equation}\label{claim}
\I(\kappa)=R[[\X]].
\end{equation}
First observe that $h^{(2)}_j\in R[\X]$, $j=1,\ldots,n$. In fact,
if $h^{(2)}_j\in R[\X,\Y]\setminus R[\X]$ for some index $j$ then
the first appearance of the $Y$-variables in $\kappa^2(Y_j)$ can
only be in a homogeneous summand of degree $>2$, contradicting the
condition $\kappa^2=\kappa$.

Assume we have shown
\begin{equation}\label{hx}
h_j^{(2)},h_j^{(3)},\ldots,h_j^{(k)}\in R[\X],\quad j=1,\ldots,n.
\end{equation}
Applying $\kappa^2=\kappa$ to the variables $X_i$ we get
\begin{equation}\label{gx}
g_i^{(2)},g_i^{(3)},\ldots,g_i^{(k)}\in R[\X],\quad i=1,\ldots,n.
\end{equation}
Therefore, (\ref{claim}) follows once it is shown that (\ref{hx})
implies
$$
h^{(k+1)}_j\in R[\X],\quad j=1,\ldots,n.
$$
Assume to the contrary that $h_j^{(k+1)}$ involves one of the
variables $Y_1,\ldots,Y_n$ for some index $j$. But then, in view
of (\ref{hx}) and (\ref{gx}), the first appearance of the
$Y$-variables in $\kappa^2(Y_j)$ can only be in a homogeneous
summand of degree higher than $k+1$, contradicting the condition
$\kappa^2=\kappa$.

The following commutative diagram is a consequence of
(\ref{claim}):
$$
\xymatrix{\hat S(P)\sqcup\hat S(Q)\ar[r]^{\subset}\ar[d]&
R[[\X,\Y]]\ar[d]_{\hat S(\alpha)}\ar[r]^{\pi\sqcup\hat S(\rho)}&
A\sqcup\hat S(Q)\ar[d]\ar[r]^\subset&R[[\X,\Y]]\ar[d]_{\hat S(\alpha)}\\
R[[\X]]\ar[r]^\subset&R[[\X,\Y]]
\ar[r]&R[[\X]]\ar[r]^\subset&R[[\X,\Y]]}
$$
where the vertical arrows represent isomorphisms, the composite of
the first two upper homomorphisms is $\pi|_{\hat S(P)}\sqcup{\bf
1}_{\hat S(Q)}$, and that of the first two lower homomorphisms is
an isomorphism. It follows that $\pi|_{\hat S(P)}:\hat S(P)\to A$
is also an isomorphism.\qed

\section{Composite automorphisms}\label{commutators}
Here we prove Theorem \ref{Main}. Let $n$ be a fixed natural
number.

Assume $\beta,\gamma\in\gac_1(R)=\Aut(R[[X]])$,
$\beta(X)=b_1X+b_2X^2+\cdots$ and $\gamma(X)=c_1X+c_2X^2+\cdots$.
We have

\begin{equation}\label{gh}
(\gamma\beta)(X)_m=\sum_{r=1}^mb_r\sum_{\substack{(j_1,\ldots,j_r)
\\ j_1+\cdots+j_r=m}}
c_{j_1}\cdots c_{j_r}.
\end{equation}

While there is no obvious compact multivariate analog of
(\ref{gh}) there is some sort of control on composite
automorphisms which will play crucial r\^ole in the sequel.

For a natural number $m$ let $\mathcal S_m$ denote the number of
ordered partitions of $m$ with nonnegative entries and having
length $n$. Arbitrary continuous $R$-algebra endomorphism $\alpha$
of $R[[\X]]$ can be represented in the following way:
\begin{align*}
\alpha(\X)=\sum_{m=1}^\infty
\alpha(\X)_m\cdot\left((X_1^{d_1}\cdots X_n^{d_n})_{
d_1+\cdots+d_n=m}\right) ^{\intercal}
\end{align*}
where $\alpha(\X)_m$ is an $n\times\mathcal S_m$ matrix over $R$
and
$$
\left((X_1^{d_1}\cdots X_n^{d_n})_{d_1+\cdots+d_n=m}\right)
$$
refers to the vector of monomials of total degree $d$ whose
components are ordered lexicographically with respect to
$X_1>\cdots>X_n$.

Next we introduce the pairings
$$
M_{n\times\mathcal S_l}(R)\times M_{n\times\mathcal S_m}(R)\to
M_{n\times\mathcal S_{lm}}(R),\quad (L,M)\mapsto L\star M,\quad
l,m\in\NN,
$$
defined by
\begin{align*}
\left(L\cdot\left((X_1^{d_1}\cdots
X_n^{d_n})_{d_1+\cdots+d_n=l}\right)^{\intercal}\right)\circ
\left(M\cdot\left((X_1^{d_1}\cdots
X_n^{d_n})_{d_1+\cdots+d_n=m}\right)^{\intercal}\right)
\left(\overline X\right)=\\
(L\star M)\cdot\left((X_1^{d_1}\cdots
X_n^{d_n})_{d_1+\cdots+d_n=lm}\right)
\end{align*}
where $\circ$ is the composition of the two endomorphisms of
$R[[\overline X]]$. The associativity of the composition yields
the associativity rule $(L\star M )\star N =L\star(M\star N)$.

Clearly, when $l=m=1$ the $\star$ operation becomes the usual
matrix product of $n\times n$ matrices.

For natural numbers $r$ and $m$ and a system of automorphisms
$\delta_1,\ldots,\delta_r\in\gac_n(R)$ the components of the
matrix $(\delta_1\cdots\delta_k)(\X)_m$ are polynomial functions
of the components of the matrices
$\delta_1(\X)_s,\ldots,\delta_k(\X)_s$, $s=1,2,\ldots$, the latter
being treated as variables in the next lemma.

\begin{lemma}\label{TOPterms}
Let $\alpha,\beta,\gamma$ be continuous $R$-algebra endomorphisms
of $R[[\overline X]]$.
\begin{itemize}
\item[(a)] The components of the matrices
$(\beta\gamma\alpha)(\X)_k$ and $(\gamma\beta)(\X)_k$ do not
depend on the components of the matrices $\alpha(\X)_m$,
$\beta(\X)_m$, or $\gamma(\X)_m$ for $k<m$.
\item[(b)] The components of the $n\times\mathcal S_m$ matrices
$$
(\beta\gamma\alpha)(\X)_m-\beta(\X)_m\star\gamma(\X)_1\star
\alpha(\X)_1-\beta(\X)_1\star\gamma(\X)_m\star
\alpha(\X)_1
$$
and
$$
(\gamma\beta)(\X)_m-\gamma(\X)_m\star
\beta(\X)_1-\gamma(\X)_1\star\beta(\X)_m
$$
do not depend on the components of the matrices $\beta(\X)_m$ or
$\gamma(\X)_m$.
\end{itemize}
\end{lemma}

Proof is straightforward.

Below for an element $\alpha\in\gic_n(R)$ we will consider the
equation $\alpha=[\gamma^{-1},\beta^{-1}]$ to be solved for
$\beta,\gamma\in\gac_n(R)$. This is equivalent to the infinite
system of equations
\begin{align*}\tag{$\mathcal E_m$}
(\beta\gamma\alpha)(\X)_m=(\gamma\beta)(\X)_m,\quad m\in\NN,
\end{align*}
for matrices $\beta(\X)_m,\gamma(\X)_m\in M_{n\times{\mathcal
S}_m}(R)$ such that $\beta(\X)_1,\gamma(\X)_1\in GL_n(R)$.

\begin{proof}[Proof of Theorem \ref{Main}] Let $b,c\in R$ be
elements such that $(b^m-1)R+(c^m-1)R=R$ for all natural numbers
$m$. For arbitrary element $\alpha\in\gic_n(R)$ we want to solve
the infinite system $(\mathcal E_m)$, $m\in\NN$.

Let $\beta(\X)_1=b\cdot{\bf Id}_n$ and $\gamma(\X)_1=c\cdot{\bf
Id}_n$. It follows from Lemma \ref{TOPterms}(a) that whatever
matrices $\beta(\X)_m$ and $\gamma(\X)_m$, $m\ge2$ we take (of
size $n\times\mathcal S_m$) the equality $(\mathcal E_1)$ is
satisfied. More generally, the same proposition implies that for
arbitrary natural number $m\ge2$ the validity of $(\mathcal E_k)$,
$k<m$ only depends on the matrices $\beta(\X)_k$ and
$\gamma(\X)_k$, $k<m$.

Now Lemma \ref{TOPterms}(b) implies that for every natural number
$m\ge2$ the equation $(\mathcal E_m)$ rewrites as
\begin{align*}
(\beta\gamma\alpha)(\X)_m-(\gamma\beta)(\X)_m=&\beta(\X)_1\star
\gamma(\X)_m\star\alpha(\X)_1+
\beta(\X)_m\star\gamma(\X)_1\star\alpha(\X)_1-\\&-\gamma(\X)_m\star
\beta(\X)_1-\gamma(\X)_1\star \beta(\X)_m+\mathcal M_m
\end{align*}
where $\mathcal M_m$ is an $n\times\mathcal S_m$ matrix which only
depends on the entries of the matrices
$\alpha(\X)_k,\beta(\X)_k,\gamma(\X)_k$, $k<m$. It is, therefore,
possible to find successively matrices $\beta(\X)_m$ and
$\gamma(\X)_m$ with the desired properties once the following is
shown -- for arbitrary natural number $m$ and arbitrary
$n\times\mathcal S_m$ matrix $\mathcal A$ there are two
$n\times\mathcal S_m$ matrices $\mathcal B$ and $\mathcal C$ such
that
$$
\mathcal A=\beta(\X)_1\star\mathcal C\star\alpha(\X)_1+\mathcal
B\star\gamma(\X)_1\star\alpha(\X)_1-\mathcal
C\star\beta(\X)_1-\gamma(\X)_1\star\mathcal B.
$$
In view of the conditions $\alpha(\X)_1={\bf Id}$,
$\beta(\X)_1=b\cdot{\bf Id}$ and $\gamma(\X)_1=c\cdot{\bf Id}$
this equation is equivalent to $ \mathcal A=\mathcal
C(b^m-b)+\mathcal B(c-c^m)$. Finally, the existence of the desired
matrices $\mathcal B$ and $\mathcal C$ follows from the condition
$$
(b^{m-1}-1)(bR)+(c^{m-1}-1)(-cR)=(b^{m-1}-1)R+(c^{m-1}-1)R=R
$$
being applied separately to every component.
\end{proof}

The proof of Theorem \ref{Main} implies the following useful fact
to be used later on:

\begin{corollary}\label{finiteapprox} Let $m$ be a natural
number. If there exist elements $b,c\in\U(R)$ such that
$(b^k-1)R+(c^k-1)R=R$ for all $k=1,\ldots,m-1$ then for arbitrary
element $\alpha\in\gic_n(R)$ there exist elements
$\beta,\gamma\in\gac_n(R)$ such that
$\alpha(\X)_k=[\beta,\gamma](\X)_k$ for $k=1,\ldots,m$.
\end{corollary}

\section{characteristic $\ge5$: the univariate case}
\label{Mushkudiani}

In this section we prove Theorem \ref{MAIN}(a). The notation
$(\mathcal E_m)$ has the same meaning as in Section
\ref{commutators}, only considered in the univariate case.

Let $\bf F$ be a field contained in $R$, $\chara{\bf F}=p\ge5$.

Choose arbitrary elements $b_2,b_3,c_2,c_3\in\bf F$ such that
$c_2-b_2\not=0$. Put $a_2=[\gamma^{-1},\beta^{-1}](X)_2$ and
$a_3=[\gamma^{-1},\beta^{-1}](X)_3$ where
$\beta,\gamma\in\gac_1(R)$ are arbitrary automorphisms such that
$\beta(X)_1=\gamma(X)=-1$ and $\beta(X)_k=b_k$, $\gamma(X)_k=c_k$,
$k=2,3$. Using (\ref{gh}) in Section \ref{commutators} one easily
sees $a_2=2(b_2-c_2)$ and $a_3=4(b_2-c_2)^2$. In particular, if
$b_2=1$, $c_2=0$ then $a_2=2$ and $a_3=4$.

\begin{lemma}\label{claimA}
Let $\alpha\in\gac_1(R)$ arbitrary element such that
$\alpha(X)=X+a_2X^2+a_3X^3+a_4X^4+\cdots$. Assume $m\ge4$ is an
even natural number, $m\not=2\mod(p)$. Then for arbitrary elements
$b_k,c_k\in R$, $k\ge4$, $k\not=m$ there exist elements
$b_m,c_m\in R$ such that $(\mathcal E_m)$ and $(\mathcal E_{m+1})$
are satisfied for the automorphisms $\beta,\gamma\in\gac_1(R)$
where $\beta(X)=-X+b_2X^2+b_3X^3+b_4X^4+\cdots$ and
$\gamma(X)=-X+c_2X^2+c_3X^3+c_4X^4+\cdots$.
\end{lemma}

\begin{proof} Using (\ref{gh}) and the
conditions on $m$ one easily checks that the system $(\mathcal
E_m,\mathcal E_{m+1})$ is independent of $b_{m+1}$ and $c_{m+1}$
(although the values of the $(m+1)$-coefficients do depend on
$b_{m+1}$ and $c_{m+1}$). Moreover, this system rewrites as the
following system of linear equations for $c_m$ and $b_m$:
\begin{align*}\tag{\bf A}
2&c_m-&2&b_m&=\mathcal X\\
\left(+2a_2-(m-2)b_2\right)&c_m+&\left(-2a_2+(m-2)c_2\right)&b_m&=\mathcal
Y
\end{align*}
where $\mathcal X$ and $\mathcal Y$ are elements of $R$ only
depending (polynomially) on $b_k,c_k$, $k<m$. See Remark
\ref{explanation} below. Now the determinant of $(\bf A)$ is
$2(m-2)(c_2-b_2)\in\U(R)$ and, therefore, the system is (uniquely)
solvable for $b_m$ and $c_m$.
\end{proof}

\begin{lemma}\label{2modp}
Assume $m\ge4$ is an even natural number, $m=2\mod(p)$. Let
$b_2=1$, $c_2=0$, and $b_3,c_3\in\bf F$ such that $b_3+c_3=0$.
Then for arbitrary elements $b_k,c_k\in R$, $k\ge4$, $k\not=m-1,m$
there exist elements $b_{m-1},c_{m-1},b_m,c_m\in R$ such that
$(\mathcal E_m)$ and $(\mathcal E_{m+1})$ are satisfied for the
automorphisms $\beta,\gamma\in\gac_1(R)$ where
$\beta(X)=-X+b_2X^2+b_3X^3+b_4X^4+\cdots$ and
$\gamma(X)=-X+c_2X^2+c_3X^3+c_4X^4+\cdots$.
\end{lemma}

\begin{proof} We start with general elements $b_2,c_2,b_3,c_3\in\bf F$.
By multiple applications of (\ref{gh}) the system $(\mathcal
E_m,\mathcal E_{m+1})$ rewrites as
\begin{align*}\tag{$\bf B$}
2(c_m-b_m)+\big(-2a_2+(m+1)b_2\big)c_{m-1}+\big(-2a_2-(m+1)c_2\big)
b_{m-1}&=\mathcal
X\\
2a_2(c_m-b_m)+\big(-3a_3+(m-4)b_3+2ma_2b_2-2a_2c_2-2b_2c_2\big)c_{m-1}
&+\\
+\big(-3a_3-(m-4)c_3+2a_2b_2-6a_2c_2+2c_2b_2\big)b_{m-1}&=\mathcal
Y
\end{align*}
where $\mathcal X,\mathcal Y\in R$ only depend (polynomially) on
$b_4,c_4,\ldots,b_{m-2},c_{m-2}$. See Remark \ref{explanation}
below. We view $(\bf B)$ as a system of linear equations with
respect to $2(c_m-b_m)$, $c_{m-1}$, and $b_{m-1}$. Since the
matrix of $(\bf B)$ is defined over $\bf F$ it is enough to show
that its rank equals 2. By multiplying the first row of the
mentioned matrix by $a_2$ and then subtracting the result from the
second row the condition on the rank becomes
$(m-3)a_2(b_2-c_2)-4b_2c_2+(m-4)(b_3+c_3)\not=0$. For our specific
choice of $b_2,c_2,b_3,c_3$ this is equivalent to
$m\not=3\mod(p)$, which is the case as $m=2\mod(p)$.
\end{proof}

\begin{remark}\label{explanation}
Here we explain how the coefficients in the systems $(\bf A)$ and
$(\bf B)$ above are computed. We only illustrate the process on
the coefficients of $c_{m-1}$ and $b_{m-1}$ in the second equation
of $(\bf B)$. The argument for all other coefficients is similar
and substantially easier.

Put $(\beta\gamma)(X)=X+d_2X^2+d_3X^3+\cdots$, the $d_j$ being
polynomial functions of $b_2,c_2,b_3,c_3,\ldots$. Then by
(\ref{gh}) we have
$$
(\beta\gamma\alpha)(X)_{m+1}=\sum_{r=1}^{m+1}a_r\sum_{\substack{(j_1,\ldots,j_r)
\\ j_1+\cdots+j_r=m+1}}d_{j_1}\cdots d_{j_r},
$$
where $d_1=1$. By Lemma \ref{TOPterms}(a) among the $d_j$ in this
equality only $d_{m-1}$, $d_m$ and $d_{m+1}$ depend on the
coefficients $c_{m-1}$ and $b_{m-1}$. Thus,
$(\beta\gamma\alpha)(X)_{m+1}=(2d_2a_2+3a_3)d_{m-1}+2a_2d_m+
d_{m+1}+\mathcal U$ where $\mathcal U$ does not depend on
$c_{m-1},b_{m-1}$.

Again (\ref{gh}), together with the assumptions that $m$ is even
and $b_1=c_1=-1$, implies  $d_2=c_2-b_2$,
$d_{m-1}=(-c_{m-1}-b_{m-1})+\mathcal U_{m-1}$,
$d_m=(m-1)b_2c_{m-1}-2c_2b_{m-1}+\mathcal U_m$ and
$d_{m+1}=\left((m-1)b_3-{{m-1}\choose2}b_2^2\right)c_{m-1}+
\left(2b_2c_2+3c_3\right)b_{m-1}+\mathcal U_{m+1}$ where $\mathcal
U_{m-1}$, $\mathcal U_m$ and $\mathcal U_{m+1}$ do not depend on
$c_{m-1}$ and $b_{m-1}$. Also, the condition $m=2\mod(p)$ implies
${{m-1}\choose2}=0\in R$.

Similar computations show
$(\gamma\beta)(X)_{m+1}=\left(3b_3+2b_2c_2\right)c_{m-1}
+(m-1)c_3b_{m-1}+\mathcal W$ where $\mathcal W$ does not depend on
$c_{m-1}$ and $b_{m-1}$. Summing up, we get the desired
coefficients in the expression
\begin{align*}
(\beta\gamma\alpha)(X)_{m+1}&-(\gamma\beta)(X)_{m+1}=\\
&\left(-3a_3+(m-4)b_3+2ma_2b_2-2a_2c_2-2b_2c_2\right)c_{m-1}+\\
&\left(-3a_3-(m-4)c_3+2a_2b_2-6a_2c_2+2c_2b_2\right)b_{m-1}+
\mathcal V+\mathcal W
\end{align*}
\end{remark}

\begin{proof}[Proof of Theorem \ref{MAIN}(a)]
Let $\alpha\in\gic_1(R)$ be an arbitrary automorphism and
$\phi\in\gac_1(R)$ be the automorphism $\phi(X)=X+2X^2+4X^3$.
Since $\chara{\bf F}\ge5$ the field $\bf F$ contains an element of
order $\ge 4$. Then by Corollary \ref{finiteapprox} (and Lemma
(\ref{TOPterms})) there exist elements
$\beta_0,\gamma_0\in\gac_1(R)$ such that
$\big(\alpha\phi^{-1}\big)(X)_k=[\beta_0,\gamma_0](X)_k$ for
$k=2,3$, equivalently
$\left([\beta_0,\gamma_0]^{-1}\circ\alpha\right)(X)_2=2$ and
$\left([\beta_0,\gamma_0]^{-1}\circ\alpha\right)(X)_3=4$. We are
done because by Lemmas \ref{claimA} and \ref{2modp} any
automorphism $\psi\in\gic_1(R)$, such that $\psi(X)_2=2$ and
$\psi(X)_3=4$, belongs to the commutator subgroup
$[\gac_1(R),\gac_1(R)]^{(1)}$.
\end{proof}

\section{characteristic $\ge5$: the multivariate case}\label{many}
In this section we assume $n\ge2$.

For a power series $g\in (X_1,\ldots,\check
X_i,\ldots,X_n)^2R[[X_1,\ldots,\check X_i,\ldots,X_n]]$ the
element $\varepsilon^g_i\in\Aut(R[\X])$, defined by
$$
\varepsilon^g_i(X_j)=
\begin{cases}
X_i+g\
\text{for}\ j=i,\\
X_j\ \text{for}\ j\not=i,
\end{cases}
$$
will be called \emph{elementary}.

\begin{lemma}\label{tamecomm} Elementary automorphisms of $R[[\X]]$
belong to $[\gac_n(R),\gac_n(R)]^{(1)}$.
\end{lemma}

\begin{proof}
It is enough to observe that
$\varepsilon^g_1\in[\gac_n(R),\gac_n(R)]^{(1)}$ for arbitrary
element $g\in (X_2,\ldots,X_n)^2R[[X_2,\ldots,X_n]]$. In fact, we
have $\beta\varepsilon^g_1\gamma=\gamma\beta$ where
$$
\beta(X_i)=
\begin{cases}
X_1+X_2\quad\text{for}\ i=1,\\
X_i\quad\text{for}\ i=2,\ldots,n
\end{cases}\quad\text{and}\quad
\gamma(X_i)=
\begin{cases}
X_i\quad\text{for}\ i=1,3,\ldots,n,\\
X_2+g\quad\text{for}\ i=2.
\end{cases}
$$
\end{proof}

\begin{lemma}\label{x1f}
Let $\frac12\in R$ and $f\in(X_2,\ldots,X_n)R[[X_2,\ldots,X_n]]$.
Then the automorphism $\alpha\in\gec_n(R)$, defined by
$$
\alpha(X_i)=
\begin{cases}
X_1+X_1f\quad\text{for}\ i=1,\\
X_i\quad\text{for}\ i=2,\ldots,n,
\end{cases}
$$
belongs to $[\gac_n(R),\gac_n(R)]^{(n)}$. Moreover, if $R$ is a
field of characteristic $\not=2$ then
$\alpha\in[\gac_n(R),\gac_n(R)]^{(1)}$.
\end{lemma}

\begin{proof} \emph{Step 1.} First consider the case when
$f\in(X_2,\ldots,X_n)^2R[[X_2,\ldots,X_n]]$. In this situation we
have $\alpha\beta\gamma=\gamma\beta$ where the elements
$\beta,\gamma\in\gac_n(R)$ are defined by
$$
\beta(X_i)=
\begin{cases}
X_1+X_1X_2\quad\text{for}\ i=1,\\
X_i\quad\text{for}\ i=2,\ldots,n
\end{cases}\quad\text{and}\quad
\gamma(X_i)=
\begin{cases}
X_i\quad\text{for}\ i=1,3,\ldots,n,\\
X_2+f+X_2f\quad\text{for}\ i=2.
\end{cases}
$$
In particular, $\alpha\in[\gac_n(R),\gac_n(R)]^{(1)}$.

\

\noindent\emph{Step 2.} Now assume $f=\xi X_2$ for some  $\xi\in
R$. Then $\alpha^2\gamma=\gamma\alpha$ where the automorphism
$\gamma\in\gac_n(R)$ is given by
$$
\gamma(X_i)=
\begin{cases}
X_i\quad\text{for}\ i=1,3,4,\ldots,n,\\
2X_2+\xi X^2_2\quad{for}\ i=2.
\end{cases}
$$
In particular, $\alpha\in[\gac_n(R),\gac_n(R)]^{(1)}$.

\

\noindent\emph{Step 3.} In the general case
$\xi_2X_2+\ldots+\xi_nX_n+g$ for some $\xi_2,\ldots,\xi_n\in R$
and $g\in(X_2,\ldots,X_n)^2R[[X_2,\ldots,X_n]]$. By iterated use
of the previous step (w.r.t. the variables $X_2,\ldots,X_n$) the
automorphism $\alpha'\in\gac_n(R)$, defined by
$$
\alpha'(X_i)=
\begin{cases}
X_1(1+\xi_2X_2)(1+\xi_3X_3)\cdots(1+\xi_n X_n)\quad\text{for}\
i=1,\\
X_i\quad\text{for}\ i=2,\ldots,n,
\end{cases}
$$
belongs to $[\gac_n(R),\gac_n(R)]^{(n-1)}$. Finally,
$\alpha=\alpha''\alpha'$ for the element $\alpha''\in\gac_n(R)$,
defined by
$$
\alpha''(X_i)=
\begin{cases}
X_1\theta\quad\text{for}\
i=1,\\
X_i\quad\text{for}\ i=2,\ldots,n,
\end{cases}
$$
where
$$
\theta=\frac{1+\xi_2X_2+\ldots+\xi_nX_n}
{(1+\xi_2X_2)(1+\xi_3X_3)\cdots(1+\xi_n
X_n)}\in1+(X_2,\ldots,X_n)^2R[[X_2,\ldots,X_n]].
$$
Now the case of rings containing 1/2 follows by Step 1.

\

\noindent\emph{Step 4.} Consider the remaining case when $R$ is a
field. By Step 1 there is no loss of generality in assuming
$f=\xi_2X_2+\cdots+x_nX_n+g$ where $\xi_2\in\U(R)$ and
$g\in(X_2,\ldots,X_n)^2R[[X_2,\ldots,X_n]]$.

If $\xi_2\not=-1$ then the exact same argument as in Step 1 works.
Therefore, we can additionally assume that $\xi_2=-1$. But then
$\alpha\beta_1\gamma_1=\gamma_1\beta_1$ for the automorphisms
$\beta_1,\gamma_1\in\gac_n(R)$ determined by
$$
\beta_1(X_i)=
\begin{cases}
X_1-X_1X_2\quad\text{for}\ i=1,\\
X_i\quad\text{for}\ i=2,\ldots,n,
\end{cases}
\ \text{and}\ \gamma_1(X_i)=
\begin{cases}
X_i\quad\text{for}\ i=1,3,4,\ldots,n,\\
X_2+(-1+X_2)f\quad\text{for}\ i=2.
\end{cases}
$$
\end{proof}

\begin{lemma}\label{h}
Assume $\ZZ_p\subset R$ and $h=(\X)^2R[[\X]]$. Then any
automorphism $\alpha\in\gac_n(R)$ of type
$$
\alpha(X_1)=
\begin{cases}
X_1+h\quad\text{for}\ i=1,\\
X_i\quad\text{for}\ i=2,\ldots,n
\end{cases}
$$
belongs to $[\gac_n(R),\gac_n(R)]^{(n+3)}$. If $R$ is a field of
characteristic $p$ then $\alpha\in[\gac_n(R),\gac_n(R)]^{(4)}$.
\end{lemma}

\begin{proof} We have $h=g+X_1f_1+X_1^2f_2+X_1^3f_3+\cdots$ for some
\begin{align*}
&g\in(X_2,\ldots,X_n)^2R[[X_2,\ldots,X_n]],\quad
f_1\in(X_2,\ldots,X_n)R[[X_2,\ldots,X_n]],\\ &f_2,f_3,\ldots\in
R[[X_2,\ldots,X_n]].
\end{align*}

First assume $g=0$. Then $\alpha=\alpha''\alpha'$ for the
automorphisms $\alpha',\alpha''\in\gac_n(R)$ defined by
$$
\alpha'(X_i)=
\begin{cases}
X_1+X_1f_1\quad\text{for}\ i=1,\\
X_i\quad\text{for}\  i=2,\ldots,n
\end{cases},\
\alpha''(X_i)=
\begin{cases}
X_1+X_1^2\theta_2+\cdots\quad\text{for}\ i=1,\\
X_i\quad\text{for}\  i=2,\ldots,n
\end{cases}
$$
where $\theta_j=(1+f_1)^{-1}f_j$, $j=2,3,\ldots$. So Theorem
\ref{MAIN}(a) and Lemma \ref{x1f} imply
$\alpha\in[\gac_n(R),\gac_n(R)]^{(n+2)}$ in general and
$\alpha\in[\gac_n(R),\gac_n(R)]^{(3)}$ for $R$ a field.

If $g\not=0$ then $\alpha=\beta\varepsilon_1^g$ where
$\beta\in\gac_n(R)$ is an automorphism of the type considered in
the preceding paragraph. So by Lemma \ref{tamecomm} we are reduced
to the previous case $g=0$.
\end{proof}

\begin{proof}[Proof of Theorem \ref{MAIN}(b,c)]
Consider an automorphism
$$
\alpha\in\gic_n(R),\quad \alpha(X_i)=X_i+h_i,\quad
h_i\in(\X)^2R[[\X]],\quad i=1,\ldots,n.
$$
We have the equality $\alpha=\alpha_1\alpha_2\cdots\alpha_n$ for
the recursively defined elements $\alpha_i\in\gac_n(R)$:
\begin{align*}
\alpha_1(X_i)&=
\begin{cases}
X_1+h_1\quad\text{for}\ i=1,\\
X_i\quad\text{for}\ i=2,\ldots,n,
\end{cases}\\
\alpha_2(X_i)&=
\begin{cases}
X_i\quad\text{for}\ i=1,3,4,\ldots,n,\\
\alpha_1^{-1}\left(X_2+h_2\right)\quad\text{for}\ i=2,
\end{cases}\\
\alpha_3(X_i)&=
\begin{cases}
X_i\quad\text{for}\ i=1,2,4,5,\ldots,n,\\
\alpha_2^{-1}\alpha_1^{-1}\left(X_3+h_3\right)\quad\text{for}\
i=3,
\end{cases}\\
\cdots&\cdots\cdots\cdots\cdots\\
\alpha_n(X_i)&=
\begin{cases}
X_i\quad\text{for}\ i=1,2,\ldots,n-1,\\
\alpha_{n-1}^{-1}\cdots\alpha_2^{-1}\alpha_1^{-1}\left(X_n+h_n\right)
\quad\text{for}\
i=n.
\end{cases}
\end{align*}
It only remains to notice that each of these automorphisms is (up
to enumeration of variables) of the type considered in Lemma
\ref{h}.
\end{proof}

\end{document}